\documentclass[12pt, reqno]{amsart} 
\usepackage{amsmath, amsthm, amscd, amsfonts, amssymb, graphicx, enumerate, color}
\usepackage[bookmarksnumbered, colorlinks, plainpages]{hyperref}
\hypersetup{colorlinks=true,linkcolor=red, anchorcolor=green, citecolor=cyan, urlcolor=red, filecolor=magenta,pdftoolbar=true}

\newtheorem{thm}{Theorem}[section]
\newtheorem{lemma}[thm]{Lemma}
\newtheorem{proposition}[thm]{Proposition}
\newtheorem{corollary}[thm]{Corollary}
\newtheorem{remark}[thm]{Remark}
\numberwithin{equation}{section}

\begin{document}
\title[Invariant mean on weakly almost periodic functionals]{ Invariant mean on weakly almost periodic functionals over  a Hopf-von Neumann algebra}

\author[R. Faal]{Ramin Faal}

\address{Department of Pure Mathematics, Ferdowsi University
of Mashhad, P. O. Box 1159, Mashhad 91775, Iran.}

\curraddr{Water, Energy, and Environmental Engineering
Research Unit, University of Oulu,
P.O. Box 4300, FI-90014 Oulu, Finland.}

\email{ramin.faal@oulu.fi}

\author[ H.R. Ebrahimi Vishki]{Hamid Reza  Ebrahimi Vishki}
\address{Department of Pure Mathematics, Ferdowsi University
of Mashhad, P. O. Box 1159, Mashhad 91775, Iran;\\ Centre of
Excellence in Analysis on Algebraic Structures (CEAAS), Ferdowsi
University of Mashhad, Mashhad, Iran.}
\email{ vishki@um.ac.ir}

\author[A. A. Khosravi]{Amin Allah Khosravi}
\address{Department of Pure Mathematics, Ferdowsi University
	of Mashhad, P. O. Box 1159, Mashhad 91775, Iran.}
\email{amin.khosravi@mail.um.ac.ir}

\subjclass[2020]{Primary: 43A60;  43A07; 46L67;  22D15.}

\keywords{Hopf-von Neumann algebra;   weakly almost periodic functional;  invariant mean.}

\begin{abstract}
Let \(M\) be a Hopf--von Neumann algebra with predual \(M_*\).
Using minimal weak*-closed left ideals, we obtain an idempotent
state \(\varphi\in\operatorname{WAP}(M_*)^*\) satisfying
\(\varphi\star \omega\star\varphi=\varphi\) for every state \(\omega\)
on \(\operatorname{WAP}(M_*)\).
This state determines the largest closed introverted subspace
\(X_\varphi\subseteq\operatorname{WAP}(M_*)\) containing the identity
and admitting a two-sided invariant mean.
We also characterize the existence of a unique two-sided invariant
mean on introverted operator systems in \(\operatorname{WAP}(M_*)\)
by a symmetry condition on their left and right state orbits.
The criterion recovers the classical results for locally compact
groups and their Fourier algebras, and gives a unified treatment
of cocommutative Hopf--von Neumann algebras, compact quantum groups,
unital \(C^*\)-Eberlein algebras, and \(\operatorname{WAP}(M(G))\)
for every locally compact group \(G\).
\end{abstract}
\maketitle

\section{Preliminaries and introduction}
\subsection{Hopf-von Neumann algebras and quantum groups}
A Hopf-von Neumann algebra is a pair $\rm (M,\Gamma )$ of a von Neumann algebra $\rm M$ and a co-multiplication $\rm \Gamma: M\to M\overline{\otimes} M$, which is a unital $*$-homomorphism satisfying $\rm (\Gamma \otimes {\rm id}) \circ \Gamma = ({\rm id}\otimes \Gamma )\circ \Gamma$, where $\overline{\otimes}$ stands for the von Neumann algebra tensor product. The main examples of a Hopf-von Neumann algebra come from locally compact groups. Indeed, for a locally compact group $\rm G$, the commutative von Neumann algebra $\rm L^\infty(G)$ and the group von Neumann algebra $\rm VN(G)$ are the main examples of Hopf-von Neumann algebras  whose co-multiplication $\rm\Gamma_a$ and $\rm \Gamma_s$ are given by
\begin{align*}
 \Gamma_a &:{\rm L^\infty (G)\to L^\infty (G)\overline{\otimes} L^\infty (G)}, \ \ {\rm\Gamma_a} (f)(g,h)=f(gh), \ \ (f\in {\rm L^\infty (G)}, g,h\in {\rm G}),\ \mbox{and}\\
\Gamma_s &:{\rm {\rm VN}(G)\to {\rm VN}(G)\overline{\otimes} {\rm VN}(G)}, \ \ {\rm \Gamma_s} (\lambda (g))=\lambda (g)\otimes \lambda (g), \ \ (g\in {\rm G}),
\end{align*}
respectively, where $\rm \lambda$ is the left regular representation of $\rm G$ on $\rm L^2(G)$. 
A somewhat less well-known example of a commutative Hopf-von Neumann  algebra is $\rm M(G)^\ast=C_0(G)^{\ast\ast}$ whose coproduct is  expertly  described in \cite[Sec. 2.1]{Dacha}. Similar to the case for $\rm L^1(G)$ and $\rm A(G)$, in  this case  also the coproduct induces the usual convolution product on the measure algebra $\rm M(G)$.
 For more information on Hopf-von Neumann algebras and locally compact quantum groups one may consult  \cite{KV2,T}.

The predual $\rm M_*$ of a Hopf-von Neumann algebra $\rm M$ enjoys  a Banach algebra structure endowed with the product $\star$ given by
\begin{align}\label{Prod}
 \langle x, \mu\star\nu\rangle =\langle \Gamma (x), \mu\otimes\nu\rangle , \ \ (x\in{\rm  M},\mu, \nu\in {\rm M_*}).
\end{align}
 The multiplication (\ref{Prod}) for the pair $({\rm L^\infty (G)},\Gamma_a)$ is the usual convolution on ${\rm L^1(G)}$, and for the pair $({\rm VN}(G),\Gamma_s)$ is the pointwise multiplication on the Fourier algebra ${\rm A(G)}$.
We also use the same notation $\star$ for the dual module operations of ${\rm M} $  as an ${\rm M} _*$-module and all of their Arens extensions. Then one can easily check  that for each $x\in {\rm M} $ and $\mu\in {\rm M} _*$, \[\mu\star x=({\rm id}\otimes \mu)\Gamma (x)\ \mbox{and}\ x\star \mu=(\mu\otimes {\rm id})\Gamma (x).\]

A very rich type of  a Hopf-von Neumann algebra is a locally compact quantum group. We provide a brief account of locally compact quantum groups whose complete description can be found in \cite{ KV2}.  A locally compact quantum group is a quadruple  $\rm (M,\Gamma,\phi,\psi)$ such that $\rm (M,\Gamma)$ is a Hopf-von Neumann algebra, $\phi$ and $\psi$ are  normal, semifinite, faithful weights, which are left and right  invariant, respectively, that is, 
\begin{align*}
\phi ((\mu\otimes {\rm id})(\Gamma x))&= \mu (1)\phi (x), \ \ (\mu\in{\rm M_*},x\in {\rm M_\phi}), \mbox{and}\\
\psi (({\rm id}\otimes \mu)(\Gamma x))&= \mu(1)\psi (x), \ \ (\mu\in{\rm M_*}, x\in {\rm M_\psi}),
\end{align*}
where  ${\rm M}_\phi := {\rm Span}\lbrace x \in M^+ : \phi (x) < \infty\rbrace$, and  ${\rm M_\psi} := {\rm Span}\lbrace x \in M^+ : \psi (x) < \infty\rbrace$.

Let $\rm G$ be a locally compact group, then the quadrapule $\mathbb{G}_a=\rm (L^\infty (G),\Gamma_a,\phi_a,\psi_a )$ is a locally compact quantum group, where $\phi_a$ and $\psi_a$ are the left and right Haar measures on $\rm G$ respectively. Furthermore, $\mathbb{G}_s=({\rm VN}(G),\Gamma_s,\phi_s)$ is a locally compact quantum group in which $\phi_s$, the Plancherel weight,   (\cite{Tak}, Deﬁnition VII.3.2) is both a left and  right Haar weight.

\subsection{Introverted subspaces and Arens products}
Let ${\rm M} $ be a Hopf-von Neumann algebra with the predual ${\rm M} _*$. Then  ${\rm M} $ is   a ${\rm M} _*$-bimodule in a natural manner. We denote   by $\star$ the multiplication of ${\rm M} _*$ and the ${\rm M} _*$-bimodule operations of ${\rm M} $.  Let ${\rm X} $ be a closed right ${\rm M} _*$-submodule of ${\rm M} $, $\varphi\in {\rm X} ^*$, $\tau\in {\rm M} _*$, and $x\in {\rm X} $. Then one can define $\tau\star \varphi :{\rm X} \to \mathbb{C}$  by $\langle\tau\star \varphi,x\rangle =\langle \varphi,x\star\tau\rangle$. When,  in addition, ${\rm X} $ is a subalgebra of ${\rm M} $ with  $1\in {\rm X} $, a functional $\varphi\in {\rm X} ^*$ with the property $\varphi (1)=\|\varphi \|=1$ is called a left invariant mean on ${\rm X} $ if
\[
\tau\star \varphi =\langle\tau,1\rangle\varphi, \ \ (\tau\in {\rm M} _*).
\]
If ${\rm X} $ has a left invariant mean, then ${\rm X} $ is called amenable.

 Let ${\rm X} $ be an ${\rm M_*} $-submodule of ${\rm M}.$  Then  ${\rm X} $ is said to be left (respectively, right) introverted, if $\omega\star x\in {\rm X} $ (respectively,  $x\star \omega\in {\rm X} $), for every $\omega\in {\rm X} ^*, x\in {\rm X} $,  where 
\[\langle \omega\star x, \mu\rangle=:\langle \omega, x\star \mu\rangle\quad\mbox{and}\quad \langle x\ast \omega, \mu\rangle=:\langle \omega, \mu\star x\rangle, \quad (\omega\in {\rm X} ^*, x\in {\rm X} , \mu\in {\rm M_*} ).\]
If ${\rm X} $ is a left (respectively, right) introverted subspace of ${\rm M}$, then ${\rm X} ^*$ equipped with the Arens type multiplication $\star$ (respectively, $\ast$) is a Banach algebra, where the operations $\ast$ and $\star$ are defined as follows:
\begin{equation}\label{Arens}
\langle \omega\star \varphi, x\rangle=:\langle \omega, \varphi\star x\rangle\quad\mbox{and}\quad \langle \omega\ast \varphi, x\rangle=:\langle \varphi, x\ast \omega\rangle, \quad (\omega,\varphi\in {\rm X} ^*, x\in {\rm X} ).
\end{equation}

\subsection{Weakly almost periodic functionals}
For a Banach algebra ${\rm A} ,$ the space of  weakly almost periodic functionals WAP$({\rm A} )$  is  defined as follows:
\[{\rm WAP}({\rm A} ):=\{f\in {\rm A} ^*:  a\mapsto a\cdot f : {\rm A} \rightarrow {\rm A} ^* \ {\rm is\ weakly\ compact}\}.
\]
It is pretty easy to verify that  WAP$({\rm A})$ is a  closed ${\rm A}$-submodule of ${\rm A} ^*$. It might be excepted that WAP$({\rm A})$  has extra structure when  ${\rm A}$ is sufficiently rich. For example, in the case where $\rm M$ is a locally compact quantum group (see \cite{KV2} for more information on locally compact quantum groups), some substantial properties of WAP$(\rm M_*)$, namely, its inclusion relation with some   other   subspaces of $\rm M$  are investigated by V. Runde \cite{Run}, that  was  continued by  Ramezanpour and Ebrahimi Vishki \cite{RV}.  In \cite{Dawe},   by a nice technique based  on  developing the theory of corepresentations on reflexive Banach spaces, for a commutative Hopf-von Neumann algebra to show that:  WAP$(\rm M_*)$ is a  $C^*$-subalgebras of $\rm M$. As a consequence, for a locally compact group $\rm G$,  $\rm WAP(M(G))$ are $\rm C^*$-algebras of $\rm M(G)^\ast$.

It is also known that, for a Hopf-von Neumann algebra $\rm M,$  ${\rm WAP}({\rm M} _*) $ is a closed  introverted operator system ${\rm M} _*$-bimodule  \cite[Proposition 5.7]{DL}.  Now, two Arens type products $\star$ and $\ast$,  as presented by \eqref{Arens}, can be defined on ${\rm WAP}({\rm M} _*)^*$; moreover,   products $\star$ and $\ast$ are equal due to \cite[Lemma 3.6]{Danon}. These two facts imply that the product $\star$ on ${\rm WAP}({\rm M} _*)^*$ is separately weak*-continuous.

Weakly almost periodic functions on a locally compact group  ${\rm G} $, denoted by ${\rm Wap}({\rm G} )$, is defined by
\[
{\rm Wap}({\rm G} )=\lbrace x\in {\rm C_b}({\rm G} ) : \lbrace L_g x : g\in {\rm G} \rbrace {\rm \ is \ relatively \ weakly \ compact \ in \ } {\rm C_b}({\rm G} ) \rbrace ,
\]
where $L_g$ denotes the left translation by $g\in {\rm G} $. 
It is known that ${\rm Wap}({\rm G} )$ is a commutative unital $C^*$-subalgebra of ${\rm C_b(G)}$ whose spectrum ${\rm G}^{\rm Wap}$ is the largest semitopological semigroup compactification of ${\rm G} $ \cite{BJM}.
From the weakly compact operators point of view, it has been shown in \cite{Ul} that $
{\rm Wap(G)=WAP(L^1(G))}.$

For a locally compact group ${\rm G}$, the commutative von Neumann algebra ${\rm L^\infty(G)}$ enjoys a natural Hopf-von Neumann algebra  structure for which
${\rm Wap(G)}$, as a right ${\rm L^1(G)}$-submodule of ${\rm L^\infty (G)}$ is amenable. Indeed, the so-called Ryll-Nardzewski fixed point theorem induces such an  invariant mean  on ${\rm Wap(G)}$;  see \cite{BJM}. Our main goal is to generalize this fact for ${\rm WAP }({\rm M} _*)$, the space of weakly almost periodic functionals on a Hopf-von Neumann algebra ${\rm M} $. The crucial key to  achieve the goal is that if ${\rm WAP }({\rm M} _*)$ admits an invariant mean, namely $\varphi$, then the left ideal generated by $\varphi$ in ${\rm WAP}({\rm M} _*)^*$,  is of dimension one.  We, therefore,  attempt to find a minimal left ideal containing a state in ${\rm WAP}({\rm M} _*)^*$  and we prove that it contains a left invariant mean. For this purpose, in Theorem \ref{main} we provide a necessary and sufficient condition for the existence of an invariant mean on ${\rm WAP}({\rm M} _*).$
We then examine it for the so-called classical Hopf-von Neumann algebras  $\rm (L^\infty(G),\Gamma_a)$ and $\rm (VN(G), \Gamma_s)$, refreshing some older results on the existence of invariant means on $\rm Wap(G)$ and $\rm Wap(\hat{G}).$ Our result also supports the existence of an invariant mean on a compact quantum group.

\section{invariant mean on ${\rm WAP}({\rm M} _*)$ }

Let   ${\rm M} $ be a Hopf-von Neumann algebra, and let ${\rm I} $ be a left ideal of ${\rm WAP}({\rm M} _*)$. Then,  the separate weak$^*$-continuity of $\star$ implies that the weak$^*$ closure of ${\rm I} $ is again a left ideal.

A state  on a  left introverted subspace $\rm X$ of $\rm M$ containing the identity $1$  of $\rm M$ is a functional  $\omega\in\rm  {\rm X} ^*$ with $\omega (1)=\|\omega\|=1$. A
 state $\omega$ is called an idempotent state if $\omega\star\omega =\omega$.
\begin{proposition}\label{minmin}
Let $\rm M$ be a Hopf-von Neumann   algebra and let $\rm X$ be an introverted subspace of ${\rm WAP (M_*)}$. Then $\rm X^*$ enjoys a minimal weak*-closed left ideal containing a  state. Furthermore, every state in this minimal left ideal is an idempotent state. 
\end{proposition}
\begin{proof}
Set
\[
\mathcal{F}=\lbrace {\rm I}\subseteq {\rm X^*}:  \mbox{{\rm I}  is a  weak*-closed  left  ideal  containing  a  state}\rbrace.
\]
It is obvious that $\mathcal{F}$ is non-empty. Let $\lbrace \rm I_\alpha\rbrace_{\alpha\in \Lambda}$ be a decreasing chain in $\mathcal{F}$ and fix a $\alpha\in\Lambda$. Then for every $\gamma\geq\alpha$,   let $\omega_\gamma\in\rm I_\gamma\subseteq\rm I_\alpha$ be a state and assume that $\omega$ be a weak*-cluster point of  the net ($\omega_\gamma$).  It follows that $\omega\in \rm I_\alpha,$ for every  $\alpha\in\Lambda,$ from which we get  $\overline{\langle \omega \rangle}^{{\rm weak}^*}\subseteq \cap_{\alpha\in\Lambda}{\rm I} _\alpha$.
Now, Zorn's lemma implies that $\mathcal{F}$ has a minimal element, say ${\rm I}$.\\

Let $\omega\in {\rm I}$ be a state which is not  an idempotent. We set 
\[
\varphi_\omega:={\rm weak^*-}\lim_n \frac{1}{n}\sum_{i=1}^n \omega^i,
\]
then $\varphi_\omega$ lies in $\rm I$ and it is  also a state with $\varphi_\omega \star \omega =\varphi_\omega =\omega\star\varphi_\omega$. 

We claim that $\omega\notin \overline{\langle \varphi_\omega\rangle}^{{\rm weak}^*}$. Indeed, if $\omega\in \overline{\langle \varphi_\omega\rangle}^{{\rm weak}^*}$, then there exists a net $(\nu_\alpha )$ such that $\nu_\alpha\star\varphi_\omega\xrightarrow{w^*} \omega$. On the other hand, since $\star$ is separately weak*-continuous we have
\begin{align*}
\nu_\alpha \star (\varphi_\omega\star\omega) \xrightarrow{{\rm weak}^*}\omega \ {\rm and} \ (\nu_\alpha \star \varphi_\omega)\star\omega \xrightarrow{{\rm weak}^*}\omega\star\omega,
\end{align*}
from which we get $\omega\star\omega=\omega$ and this is a contradiction. We thus arrive at $\overline{\langle \varphi_\omega\rangle}^{w^*}$ as a proper left ideal in ${\rm I} $.
\end{proof}

In the next result we need to benefit from Banach-Dieudonn\'{e} theorem stating that a convex subset $C$ of a dual Banach space ${\rm X} ^*$ is weak*-closed if and only if $C\cap n {\rm Ball}({\rm X} ^*)$ is weak*-closed for every $n\in\mathbb{N}; $ \cite[Theorem 3.92]{Fab}.
\begin{proposition}\label{cl}
Let $\rm M$ be a Hopf-von Neumann  algebra and let $\rm X$ be an introverted subspace of ${\rm WAP(M_*)}$. Then for every idempotent state $\omega$ in $\rm X^*$, the left ideal  $\langle \omega \rangle$ is weak*-closed.
\end{proposition}
\begin{proof}
By Banach-Dieudonn\'{e} theorem, as $\langle \omega\rangle=\rm X^*\star\omega$ is a subspace, we just need to show that $\langle \omega \rangle\cap {\rm Ball}({\rm X}^*)$ is weak*-closed.
Let $(\nu_\alpha\star \omega)$ be an arbitrary net in $\langle \omega \rangle\cap {\rm Ball}({\rm X}^*),$ weak*-converging to $\nu.$
We will show that $\nu\in\langle \omega\rangle$. Since   $\|\nu_\alpha\star \omega \|\leq 1$,  for every $\alpha$, we get
\begin{align*}
\|{\nu_\alpha}|_{\omega\star {\rm X}}\| &=\sup_{\|\omega\star x \|\leq 1}|\langle \nu_\alpha , \omega\star x\rangle | 
=\sup_{\|\omega\star x \|\leq 1}|\langle \nu_\alpha\star \omega , \omega\star x\rangle |  \\
&\leq \sup_{\|\omega\star x \|\leq 1} \|  \nu_\alpha\star \omega \| \| \omega\star x \|  \leq 1,
\end{align*}
in which, $\omega\star {\rm X} =\{ \omega\star x : x\in {\rm X}\}$. An standard argument based on the Hahn-Banach theorem now implies that, for each $\alpha$, we can extend  ${\nu_\alpha}|_{\omega\star {\rm X}}$ to an element $\tilde{\nu}_\alpha\in {\rm Ball}({\rm X}^*)$, and so   $\nu_\alpha\star\omega =\tilde{\nu}_\alpha\star\omega$, for every $\alpha.$  
Considering $\tilde{\nu}\in {\rm Ball}({\rm X}^*)$ as a weak*-cluster point of  the net $(\tilde{\nu}_{\alpha})$, then from the separate weak*-continuity of $\star$ we arrive at 
$\nu=\tilde{\nu}\star\omega\in \langle \omega\rangle$, as claimed.
\end{proof}

To arrive at the main result, we need the following key lemma. 
\begin{lemma}\label{P}
Let $\rm M$ be a Hopf-von Neumann algebra and $\rm X$ be an introverted operator system  of $\rm Wap(M_*)$ containing the identity of  $\rm M.$ Then there exists an idempotent state $\varphi$ in  $\rm X^*$ such that   $\varphi =\varphi\star \omega\star \varphi$ for every state $\omega\in\rm X^*.$
\end{lemma}
\begin{proof}
 Let $\rm I$ be  the minimal  weak*-closed left ideal  of $\rm X^*$ containing an idempotent state $\varphi$; see Proposition \ref{minmin}.  Then for every state  $\omega\in\rm X^*,$  it follows from   Proposition \ref{minmin} that $\omega\star\varphi$ is an idempotent state and the minimality of $\rm I$ together with Proposition  \ref{cl} imply that  $\rm I=\langle\omega\star\varphi\rangle.$ Now   $\varphi\in\rm I=\langle\omega\star\varphi\rangle$ and there exists  a $\nu\in\rm  X^*$ such that $\nu\star\omega\star\varphi=\varphi$.

\[
\varphi\star\omega \star\varphi = (\nu\star\omega\star\varphi) \star (\omega\star\varphi) =\nu\star (\omega\star\varphi\star \omega\star\varphi) = \nu\star\omega\star\varphi=\varphi,
\]
as claimed.
\end{proof}

We use Lemma \ref{P} to give the following result which also implies  the existence of a unique invariant mean on a compact quantum group {\rm {\cite[Theorem 2.4]{Van}}}, cocommutative Hopf-von Neumann algebras, and Eberlin $C^*-$algebras independently.
\begin{corollary}\label{SI}
If $\rm M$ is a Hopf-von Neumann algebra,  and for every two idempotent states $\omega$ and $\omega'$  in  ${\rm WAP(M_*)^*}$ the identity $\omega\star \omega' =\omega'$ holds if and only if $\omega'\star\omega =\omega'$, then ${\rm WAP(M_*)}$ admits a unique invariant mean. 
\end{corollary}
\begin{proof}
Lemma \ref{P} yields that there exists an idempotent state $\varphi$ such that for each state $\omega$ we have $\varphi\star\omega\star\varphi =\varphi$, and this implies that  $\varphi\star\omega = \varphi$, that is,   $\varphi$ is a right invariant mean.  Similarly, one can find a left invariant mean $\varphi'$ and thus $\psi=\varphi'\star\varphi$ is an invariant mean, whose uniqueness is clear.
 \end{proof}

\begin{remark}
\begin{enumerate}
\item 
Let  $\mathbb{G}$ be a  compact quantum group, then ${\rm WAP(L^1(\mathbb{G}))=C(\mathbb{G})}$, and \cite[Proposition 2.8]{Salsa} implies that the desired condition holds for every two idempotent states $\omega$ and $\omega'$  in  ${\rm WAP(L^1(\mathbb{G}))^*=M(\mathbb{G})}$. 

    \item 
$C^*$-Eberlein algebras were introduced in
\cite[Definition 3.6]{DaDa}. These are $C^*$-algebras
equipped with a coproduct and a contractive Hilbert-space
generator whose coefficients are norm dense. Let the quadruple $(A,\Phi,V,H)$
is a \(C^*\)-Eberlein algebra, and assume for now that \(A\) is unital.
By the definition presented in \cite{DaDa}, the map
\[
    \pi:A^*\longrightarrow B(H),
    \qquad
    \pi(\mu)=(\mu\otimes\mathrm{id})(V)
\]
is a contractive homomorphism, and the coefficients
\[
   \{ (\mathrm{id}\otimes\omega)(V),
    \quad \omega\in B(H)_*\},
\]
are norm dense in \(A\). The density of the coefficients implies that
\(\pi\) is injective. In fact, if \(\pi(\mu)=0\), then
    $\mu\bigl((\mathrm{id}\otimes\omega)(V)\bigr)
    =\omega\bigl(\pi(\mu)\bigr)=0$
for every \(\omega\in B(H)_*\), and hence \(\mu=0\).
Now let \(e,f\in A^*\) be idempotent states, and put
    $P=\pi(e),
    Q=\pi(f).$
Since \(\pi\) is contractive and multiplicative, we get
\[
    P^2=P,
    \qquad
    Q^2=Q,
    \qquad
    \|P\|,\|Q\|\leq 1.
\]
Every contractive idempotent on a Hilbert space is an orthogonal
projection.
Thus \(P=P^*\) and \(Q=Q^*\), so
\[
\begin{aligned}
    e\star f=f
    &\Longleftrightarrow PQ=Q\\
    &\Longleftrightarrow (PQ)^*=Q^*\\
    &\Longleftrightarrow QP=Q\\
    &\Longleftrightarrow f\star e=f.
\end{aligned}
\]
Therefore, $e\star f=f$ if and only if $f\star e=f$.
This is precisely the condition in Corollary \ref{SI}.

To place this example within the von Neumann algebraic framework, take
$    M=A^{**}, X=A\subseteq M,$
and extend \(\Phi\) normally to \(M\).
The Hilbert space coefficients show that
    $A\subseteq
    \operatorname{WAP}(M_*)
    =
    \operatorname{WAP}(A^*)$.
Moreover, the module actions preserve the coefficients. Consequently,
\(X=A\) is an introverted operator system.
Thus the same argument as in Corollary \ref{SI}, with \(X=A\), yields a unique
invariant mean.

In particular, if \(\varphi\) be as in Lemma \ref{P} then
\[
    \varphi\star \omega\star\varphi=\varphi
    \qquad (\omega\in S(A)),
\]
and, setting \(P=\pi(\varphi)\) and \(T=\pi(\omega)\), we obtain $PTP=P.$
For \(\xi\in PH\),
\[
    PT\xi=\xi,
    \qquad
    \|T\xi\|\leq\|\xi\|.
\]
Since \(P\) is an orthogonal projection, this forces \(T\xi=\xi\).
Applying the same argument to \(T^*\), we obtain $TP=PT=P.$
The injectivity of \(\pi\) now gives
        $\omega\star\varphi=\varphi\star \omega=\varphi.$

Moreover, concerning the Eberlein $\rm C^*$-algebra $\rm E(\mathbb G)$, 
Theorem 3.28 and Remark 3.30 of \cite{DaDa} imply that if we work with $\rm C_0(\mathbb{G})$ instead of $\rm C_0^u(\mathbb{G})$ and construct $\rm E(\mathbb{G})$ using unitary corepresentations of $\rm C_0(\mathbb{G})$ as a $\rm C^*$-subalgebra of $\rm C_b(\mathbb{G})$, then it would lie in ${\rm WAP(M(\mathbb{G}))}$. Now, the fact that  $\rm E(\mathbb{G})\subseteq {\rm WAP(M(\mathbb{G}))\cap C_b(\mathbb{G})}$ together with the fact that  $\rm L^1(\mathbb{G})$ can be embedded in $\rm C_0(\mathbb{G})^*$ as a closed two sided ideal, implies that $\rm E(\mathbb{G})\subseteq {\rm WAP}(L^1(\mathbb{G}))$. See {\rm (\cite[Theorem 5.3]{DaDa})} for admitting a unique invariant mean of $\rm E(\mathbb G)$.
\end{enumerate}
\end{remark}

In the following, we characterize the largest closed introverted operator system containing $1$ with the virtue of $\varphi$ obtained in Lemma \ref{P}.

\begin{thm}\label{largest}
Let ${\rm M}$ be a Hopf-von Neumann algebra, and let
$\varphi\in {\rm WAP(M_*)}^*$ be the idempotent state
provided by Lemma \ref{P}, applied to ${\rm WAP(M_*)}$.
Let $P$ be the set of all states in ${\rm WAP(M_*)}^*$,
and define
\[
X_\varphi
=\bigcap_{p,q\in P}
\ker\!\left(p\star\varphi\star q-\varphi\right).
\]
Then $X_\varphi$ is a closed introverted operator system
containing $1$, and $\varphi|_{X_\varphi}$ is its unique
two-sided invariant mean. Moreover, $X_\varphi$ is the largest closed introverted
subspace of ${\rm WAP(M_*)}$ containing $1$ and admitting
a two-sided invariant mean.
\end{thm}

\begin{proof}
By its definition, $X_\varphi$ is a closed linear subspace
of ${\rm WAP(M_*)}$. Since $p\star\varphi\star q$ and
$\varphi$ are states, their difference annihilates $1$
and is a self-adjoint functional. Hence $1\in X_\varphi$
and $X_\varphi$ is self-adjoint. It is also evident that $P\star X_\phi, X_\phi\star P\subseteq X_\phi$. Since states linearly span ${\rm WAP(M_*)}^*$, the same
holds for every $\rho\in{\rm WAP(M_*)}^*$.
Every functional on $X_\varphi$ has a Hahn--Banach
extension to ${\rm WAP(M_*)}^*$, so $X_\varphi$ is
introverted.

Put $m_\varphi=\varphi|_{X_\varphi}$.
Since $1\in X_\varphi$, we get
$\|m_\varphi\|=m_\varphi(1)=1.$ As $\varphi$ is an idempotent state, we get 
\[
\langle p\star\varphi,x\rangle =\langle p\star\varphi\star\varphi,x\rangle
=\langle\varphi,x\rangle
=\langle\varphi\star\varphi\star p,x\rangle = \langle\varphi\star p,x\rangle
\qquad(p\in P,\ x\in X_\varphi).
\]
by the definition of $X_\phi$.
By linearity and Hahn--Banach extension, these identities
imply
\[
\nu\star m_\varphi
=\langle\nu,1\rangle m_\varphi
=m_\varphi\star\nu
\qquad(\nu\in X_\varphi^*).
\]
Therefore $m_\varphi$ is a two-sided invariant mean.
If $n$ is another two-sided invariant mean, then
$m_\varphi\star n=n$
 and 
$m_\varphi\star n=m_\varphi,$
so $n=m_\varphi$.

Finally, let ${\rm Y}$ be a closed introverted subspace
of ${\rm WAP(M_*)}$ containing $1$, and let $m$ be a
two-sided invariant mean on ${\rm Y}$.
The restriction map
\[
r:{\rm WAP(M_*)}^*\longrightarrow {\rm Y}^*,
\qquad r(\omega)=\omega|_{\rm Y},
\]
is a surjective homomorphism by \eqref{Arens}.
Choose a norm-preserving Hahn--Banach extension
$\mu\in P$ of $m$, and put $s=r(\varphi)$.
Lemma \ref{P} yields
\[
s\star m\star s
=r(\varphi\star\mu\star\varphi)
=r(\varphi)=s.
\]
Since $s$ is a state, the two-sided invariance of $m$
also gives $s\star m\star s=m$. Consequently,
$r(\varphi)=m$.
For every $p,q\in P$ and $y\in{\rm Y}$, we therefore have
\[
\begin{aligned}
\langle p\star\varphi\star q,y\rangle
&=\langle r(p)\star m\star r(q),y\rangle\\
&=\langle m,y\rangle
=\langle\varphi,y\rangle.
\end{aligned}
\]
Hence
\[
{\rm Y}\subseteq
\bigcap_{p,q\in P}
\ker\!\left(p\star\varphi\star q-\varphi\right)
=X_\varphi,
\]
which proves maximality.
\end{proof}

We are ready to state and prove our result providing an equivalence condition to the existence and uniqueness of an invariant mean on $\rm WAP(M_*),$ for a  Hopf-von Neumann algebra $\rm M$.
\begin{thm}\label{main}
Let $\rm M$ be a Hopf-von Neumann algebra  and let $\rm X$ be an introverted operator system in  $\rm WAP(M_*)$ containing the identity $\rm 1$  of $\rm M.$ Let $P$ be the set of all states $\nu$ in $\rm X^*$ with $\nu(1)=1$.  Then the following statements are equivalent.
\begin{enumerate}[\hspace{1em}\rm (a)]

\item[$(a)$]
$\rm X^*$ has an invariant mean (necessarily unique).
\item[$(b)$]
${\rm P}\star x$ contains non-zero multiplication of identity if and only if $x\star {\rm P}$ contains (the same) non-zero multiplication of identity, for every $x\in\rm X$.
\end{enumerate}
\end{thm}
\begin{proof}
Let $\rm X^*$ has an invariant mean, namely $\varphi$. Then $\varphi\in \rm P$ and for every $\tau\in {\rm M}_*$,
\[
\langle \varphi\star x,\tau\rangle =\langle  \varphi\star\tau,x\rangle = \langle\tau,1\rangle \langle \varphi , x\rangle.
\]
Thus $\varphi\star x= \langle \varphi, x\rangle 1$. When $\langle \varphi,x\rangle\neq 0$, we have $\alpha 1\in {\rm P}\star x$ for some non-zero scalar $\alpha$. In the case when $\langle \varphi, x\rangle = 0$, then ${\rm P}\star x$ and $x\star {\rm P}$ do not contain any non-zero  scalar multiplication of identity. Indeed, if $\alpha 1=x\star \mu$ for some non-zero scalar $\alpha$ and $\mu\in \rm P$, then
\[
 \alpha 1 = \alpha \langle \varphi, 1\rangle = \langle \varphi, x\star \mu\rangle =\langle \mu, \varphi\star x\rangle =0.
\]
  Therefor, The  statement $(a)$ follows $(b)$.
  
   To prove the converse, let  $(b)$ holds for every $x\in {\rm X} .$ 
Suppose that  $\varphi$  is  an idempotent state in  $\rm P$  satisfying $\varphi=\varphi\star\omega\star\varphi$, for every state $\omega\in\rm X^*,$ as appeared in Lemma \ref{P}. If  for every $x\in\rm X,$ $\varphi\star x\star\varphi=\varphi\star x,$ 
then the latter identity  implies  $\omega \star \varphi = \varphi $, for every state  $\omega\in {\rm X} ^*$, and this confirms that $\varphi$ is a left invariant mean.
On the contrary, let $z =\varphi\star x\star\varphi -\varphi\star x \neq 0$ for some $x\in\rm X$, and set 
\[
{\rm I}_z := \{ \omega\in{\rm P} : \omega\star z =0\}.
\]
Considering $\rm P$ as weak*-closed convex semigroup, ${\rm I}_z$ is a  left ideal excluding $\varphi$ which is weak*-closed. We show that ${\rm I}_z$ contains a state.

 First note that, $1\notin{\rm P}\star z$. Indeed, if $1\in{\rm P}\star z$, then for some $\omega\in\rm P$ we get $1 = \omega\star z$. Consequently, 
\[
1 = 1 \star\varphi=  \omega\star z \star \varphi=0.
\]
Thus, $1\notin z\star\rm P$ by the assumption $(b)$. As $z\in \rm WAP(M_*)$, $z\star\rm P$ is norm closed. In fact, let $(z\star \omega_n)$ be a sequence in $z\star\rm P$ which tends to $z_0$, and let $(\omega_{n_k})$ converges weak* to $\omega\in \rm P$. then, considering  $z\in \rm WAP(M_*)$,  $(z\star\omega_{n_k})$ tends weakly to $z\star\omega$. This follows that $z_0=z\star\omega\in z\star\rm P$. 

Since $z\star \rm P$ is convex, norm closed, and exclude $1$, there exists   a functional $\tau\in X^*$ with $\langle \tau,1\rangle =1$ and $\tau|_{z\star {\rm P}}= 0$. The restriction of a state on $M$ to $X$ induce a state on $X$, and this follows that each functional on $X$ is a linear span of some elements in $\rm P$. Thus we get $\tau\star z =0$ and $\tau\star\varphi\neq 0$. As $\varphi \star \tau\star\varphi =\langle\tau, 1\rangle\varphi$, the equality $\tau\star\varphi \star \tau\star\varphi =\langle\tau, 1\rangle\tau \star \varphi$ implies that $\tau\star \varphi$ is an idempotent state. Furthermore, the identity $\tau\star \varphi\star z = \tau\star z =0$ deduces that $\tau\star \varphi\in {\rm I} _z$.
 Let $I_\varphi$ be the left ideal in $\rm P$ generated by $\varphi$. Since $\varphi$ is idempotent, $I_\varphi$ is weak*-closed and contains $\tau\star\varphi\star\varphi=\tau\star\varphi$. Thus ${\rm I} _z\cap \langle\varphi\rangle$ is a weak*-closed proper left ideal (excluding $\varphi$) containing a state (including $\tau\star\varphi$),  which contradicts the minimality of $I_\varphi$ . Therefore   $\varphi\star x\star\varphi=\varphi\star x,$ for every $x\in\rm X,$ and as it was stated  $\varphi$ is a left invariant mean. In the same fashion, one can find a right invariant mean $\varphi'$ and thus $\psi=\varphi\star\varphi'$ is an invariant mean, whose uniqueness is clear. 
\end{proof}

It is worthwhile mentioning that the condition $\rm (b)$ of Theorem \ref{main} holds in the case where the Hopf-von Neumann algebra $\rm M$ is cocommutative. Indeed, in this case $\rm M_*$ is commutative,  and so  for every $x\in\rm WAP(M_*)$ and $\tau\in M_*$, $x\star \tau=\tau\star x$. Now for every  introverted operator system in  $\rm WAP(M_*)$ containing the identity of $\rm M,$  an standard argument  reveals that $x\star \omega=\omega\star x,$ for every $x\in\rm X$ and $\omega\in X^*$. We arrive at the following corollary.
\begin{corollary}\label{cocom}
For every cocommutative  Hopf-von Neumann algebra  $\rm M,$ $\rm WAP(M_*)$ admits a unique invariant mean.
\end{corollary}
  Applying  Corollary \ref{cocom}  for the  so-called classical cocommutative Hopf-von Neumann algebras $(\rm VN(G),\Gamma_s)$, on a locally compact group $\rm G$,  we get  the following known result.
\begin{corollary}\label{classic}{\rm({\cite[Corollary 5.7]{La}})}. Let $\rm G$ be a locally compact group, then
  ${\rm Wap(\widehat{G})}={\rm WAP(A(G))}$ admits a unique invariant mean.

  The question of existence and uniqueness of an invariant mean on ${\rm WAP(M(G))}$, on a locally compact group $\rm G$,  raised by  Dales,  Lau and Strauss  \cite[Sec. 10]{seme}, however, the same existence problem has already formulated by  Daws  \cite[pp. 222 and 233]{Dacha}. The problem fully answered in affirmative for every locally compact group ${\rm G} $ by Neufang \cite[Theorem 2.3]{Neuf}. However, in the case where $\rm G$ is abelian, the result can be also derived by applying Corollary \ref{cocom} for the cocommutative Hofp von Neumann algebra $(\rm M(G)^*, \Gamma)$, whose coproduct $\Gamma$  was described in \cite[pp. 288-289]{Dawe}.
\end{corollary}

\begin{corollary}
	${\rm Wap(G)}$ admits a unique invariant mean for each locally compact group $G$.
\end{corollary}
\begin{proof}
Set
Let \(f\in WAP(G)\) and $P=S(WAP(G))$. consider
\[
K_L(f)=\overline{\operatorname{co}\{L_sf:s\in G\}}^{\|\cdot\|},
\qquad
K_R(f)=\overline{\operatorname{co}\{R_sf:s\in G\}}^{\|\cdot\|}.
\]
Both sets are weakly compact. The corresponding translations act on them as weakly continuous affine isometries. The Ryll–Nardzewski fixed-point theorem therefore gives a fixed point in each set. A function fixed by every left translation, or every right translation, is constant. Thus both sets contain constants. 
Point evaluations at elements of \(G\) are states on \(WAP(G)\). Since \(P\star f\) and \(f\star P\) are norm closed, it follows that each contains a constant.
Now suppose for some non-zero $\alpha$ and $\nu\in P$, we have
$f\star\nu=\alpha1$. Choose \(\mu\in P\) and \(\beta\in\mathbb C\) such that $\mu\star f=\beta1.$ 
The commuting left and right module actions give
\[
\alpha1
=\mu\star(f\star\nu)
=(\mu\star f)\star\nu
=\beta1.
\]
Hence \(\beta=\alpha\ne0\). The converse is identical. This proves condition (b) of Theorem \ref{main}.
\end{proof}
\begin{corollary}

\(\operatorname{WAP}(M(G))\) admits a unique invariant mean for every locally compact group $G$.
\end{corollary}

\begin{proof}
Fix \(f\in \operatorname{WAP}(M(G))\), and choose \(a\in L^1(G)_+\) with
\(\|a\|_1=1\).
The smoothed functionals \(a\star f\) and \(f\star a\) are
represented, respectively, by
\[
 s\longmapsto f(\delta_s*a),
 \qquad
 s\longmapsto f(a*\delta_s).
\]
These functions belong to \(\operatorname{WAP}(G)\). Their
Continuity follows from norm continuity of translations on
\(L^1(G)\), and weak compactness of their translation orbits
is inherited from \(f\).

Ryll--Nardzewski gives constants in the norm-closed convex
right translation hull of \(a\star f\) and left translation
hull of \(f\star a\).
The translates act on \(f\) through the probability measures
\(\delta_s*a\) and \(a*\delta_s\), respectively, whose restrictions
to \(\operatorname{WAP}(M(G))\) are states.
Since \(P\star f\) and \(f\star P\) are convex and norm closed,
there exist \(\mu,\nu\in P\) and
\(\alpha,\beta\in\mathbb C\) such that $\mu\star f=\alpha 1$ and
 $f\star\nu=\beta 1.$
Moreover,

\[
 \beta 1
 =\mu\star(f\star\nu)
 =(\mu\star f)\star\nu
 =\alpha 1.
\]
Thus $\operatorname{WAP}(M(G))$ satisfied the condition Theorem \ref{main} $(b)$.
\end{proof}

\begin{remark}\label{isometric-spaces}
We use ${\rm WAP}_{\rm iso}(M)$ in the sense of
\cite[Definitions 2.2, 2.7, 2.8 and 2.11]{Yul}.
A left covariant pair on a reflexive Banach space $V$ consists
of contractive linear maps
\[
\pi:M_*\longrightarrow B(V),
\qquad
\xi:M_*\longrightarrow V,
\]
where $\pi$ is a homomorphism and
\[
\xi(\mu\star\nu)=\pi(\mu)\xi(\nu)
\qquad(\mu,\nu\in M_*).
\]
Its coefficients are the elements $x\in M$ satisfying
\[
\langle x,\mu\rangle=\langle\xi(\mu),\eta\rangle
\qquad(\mu\in M_*)
\]
for some $\eta\in V^*$.
The space ${\rm WAP}_{\rm iso}(M)$ consists of these coefficients
when, in addition, $\pi$ is $\ell^2$-contractive and has a
surjective isometric generator in some faithful normal spatial
realization of $M$, as in the cited definitions.
Here $\ell^2$-contractivity means
\[
\bigl\|[\pi(\mu_{ij})]\bigr\|_{B(\ell^2_n(V))}
\leq \bigl\|[\mu_{ij}]\bigr\|_{M_n(M_*)}
\qquad(n\in\mathbb N),
\]
where $\ell^2_n(V)$ is equipped with the norm
\[
\|(v_j)_{j=1}^n\|
=\left(\sum_{j=1}^n\|v_j\|^2\right)^{1/2}.
\]
In particular,
${\rm WAP}_{\rm iso}(M)\subseteq{\rm WAP}(M_*)$.
Suppose that in addition
\[
\Gamma(z)=1\otimes z\quad\Longrightarrow\quad z\in\mathbb C1
\qquad(z\in M).
\]
This is the right nontriviality condition of
\cite[Definition 3.2]{Yul}, and it holds for locally compact
quantum groups. Let $P_*$ denote the normal states of $M$, and write
\[
K_L(x)=\overline{P_*\star x}^{\|\cdot\|},
\qquad
K_R(x)=\overline{x\star P_*}^{\|\cdot\|}.
\]
The normal state orbits are convex, so their norm and weak
closures coincide. By
\cite[Propositions 3.1 and 3.4]{Yul}, both $K_L(x)$ and $K_R(x)$
contain a scalar multiple of $1$ whenever
$x\in{\rm WAP}_{\rm iso}(M)$.
These scalar multiples are unique and equal. Put 
\[
Y_{\rm iso}
=\overline{{\rm WAP}_{\rm iso}(M)
+{\rm WAP}_{\rm iso}(M)^*}^{\|\cdot\|}.
\]
Then it can be checked that $Y_{\rm iso}$ is a closed self-adjoint introverted $M_*$-bimodule in
${\rm WAP}(M_*)$ containing $1$. Moreover, condition {\rm (b)} of Theorem \ref{main} for $Y_{\rm iso}$ can be deduced from
\cite[Propositions 3.1 and 3.4]{Yul}.

\end{remark}

\begin{remark}\label{strict-examples}
The inclusion
${\rm WAP}_{\rm iso}(M)\subseteq X_\varphi$
can be strict.
Let
\[
S=(\mathbb N_0,\max),\qquad
M=\ell^\infty(S),\qquad
\Gamma(f)(m,n)=f(\max\{m,n\}).
\]
This is a cocommutative Hopf-von Neumann algebra,
and $M_*=\ell^1(S)$ has identity $\delta_0$.
Let $c$ be the space of convergent sequences.
Then ${\rm WAP}(M_*)=c.$
Indeed, for
\[
T_nf(k)=f(\max\{k,n\}),
\]
if $f$ converges to $L$, then $T_nf\to L1$ in norm,
so the translation orbit is relatively compact.
Conversely, if $f$ is weakly almost periodic, a
subnet of $(T_nf)$ converges weakly to a constant
$L1$. Since $T_nf-f\in c_0$ and $c_0$ is weakly
closed, $L1-f\in c_0$, proving $f\in c$.

The functional
\[
\varphi(f)=\lim_{n\to\infty}f(n)
\qquad(f\in c)
\]
is an invariant mean. Its uniqueness follows from
$T_nf\to\varphi(f)1$ in norm. Hence Lemma
\ref{largest} gives $X_\varphi=c$.

We claim that ${\rm WAP}_{\rm iso}(M)=\mathbb C1$.
Let $(\pi,\xi)$ be a defining covariant pair with
isometric generator $U$ in a faithful normal
realization $\sigma:M\to B(H)$.
For $p_n=1_{\{n\}}$, choose a unit vector
$h_n\in\sigma(p_n)H$.
The generator identity gives
\[
U(h_n\otimes v)=h_n\otimes\pi(\delta_n)v
\qquad(v\in V).
\]
Thus $\pi(\delta_n)$ is an isometry.
Since $\delta_n\star\delta_n=\delta_n$, it is also
idempotent, and therefore $\pi(\delta_n)=I_V$.
By linearity and continuity,
\[
\pi(\mu)=\langle\mu,1\rangle I_V
\qquad(\mu\in\ell^1(S)).
\]
Covariance and the identity $\delta_0$ now imply
\[
\xi(\mu)=\pi(\mu)\xi(\delta_0)
=\langle\mu,1\rangle\xi(\delta_0).
\]
Every coefficient is consequently constant.
The reverse inclusion follows from the trivial
representation.
We have therefore obtained
\[
\overline{{\rm WAP}_{\rm iso}(M)}^{\|\cdot\|}
=
{\rm WAP}_{\rm iso}(M)
=
\mathbb C1
\subsetneq
X_\varphi
=
{\rm WAP}(M_*)
=
c.
\]
\end{remark}

\begin{remark}
\begin{enumerate}[\hspace{1em}\rm (1)]

\item  In the realm of locally compact groups (see \cite{KV2}),   Daws \cite{Danon} introduced ${\rm Wap}(\mathbb{G})$, for a locally compact quantum group $\mathbb{G}$, as the biggest  $C^*$-algebra inside ${\rm WAP(L^1(\mathbb{G}))}$. More precisely, he defined  ${\rm Wap}(\mathbb{G})$ by
\[
{\rm Wap}(\mathbb{G})=\lbrace x\in {\rm WAP(L^1(\mathbb{G}))} : x^*x,xx^*\in {\rm WAP(L^1(\mathbb{G}))}\rbrace,
\]
and showed that ${\rm Wap}(\mathbb{G})$ is an introverted $\rm L^1(\mathbb{G})$-bimodule $C^*$-subalgebra of $\rm L^\infty (\mathbb{G})$. Consequently, he raised 
 a question  (see  \cite[Sec. 7, second question]{Danon}) concerning the existence of an invariant mean on ${\rm Wap}(\mathbb{G})$.

The existing examples supports the conjecture that  ${\rm Wap}(\mathbb{G})$ satisfies the condition $\rm (b)$ of Theorem \ref{main}.

\item Concerning the conclusion of Theorem \ref{main},  we leave the following question undecided:

\textit{Is there a Hopf-von Neumann algebra $\rm M$ for which ${\rm WAP(M_*)}$
does not enjoy an invariant mean?}
\end{enumerate}

\end{remark}

\begin{remark}
 It is worthwhile mentioning that all of our arguments and so the  results, particularly the main one (Theorem \ref{main}), hold for an arbitrary  Lau algebra $\rm A$  instead of the predual $\rm M_*$ of a Hopf-von Neumann algebra. 
We recall that  a Lau algebra (same as F-algebra in \cite{La1})  is a  Banach algebra
$\rm A$ such that $\rm A^*$  is a von Neumann algebra whose unit  lies in the spectrum of $\rm A$. 
\end{remark}
\section*{acknowledgment} 
\thanks{This work was mainly carried out while the author was
affiliated with Ferdowsi University of Mashhad, Mashhad, Iran,
and was completed at the University of Oulu, Oulu, Finland.}

\end{document}